\documentclass[12pt]{amsart}
\usepackage{amsmath, amssymb}
\usepackage{amsthm,latexsym}

\newtheorem{theorem}{Theorem}
\newtheorem{corol}{Corollary}
\newtheorem{lemma}{Lemma}
\newtheorem{propos}{Proposition}

\author{Vitaly V. Balashchenko}

\title{Invariant $f$-structures in the
generalized Hermitian geometry}

\date{}


\begin{document}

\maketitle
\begin{abstract}

We collect the recent results on invariant $f$-structures in the
generalized Hermitian geometry. Here the canonical $f$-structures
on homogeneous $k$-symmetric spaces play a remarkable role.
Specifically, these structures provide a wealth of invariant
examples for the classes of nearly K\"ahler $f$-structures,
Hermitian $f$-structures and some others. Finally, we consider all
invariant $f$-structures on the complex flag manifold
$SU(3)/T_{max}$ and describe them in the sense of generalized
Hermitian geometry. In particular, it presents first invariant
examples of Killing $f$-structures.
\end{abstract}


\tableofcontents 

\section{\bf Introduction}

Invariant structures on homogeneous manifolds are traditionally
one of the most important objects in differential geometry,
specifically, in Hermitian geometry. Some remarkable classes of
almost Hermitian structures such as K\"ahler, nearly K\"ahler,
Hermitian structures etc. are well known and intensively used in
geometry and a number of applications. In particular, a special
role is played by a significant class of invariant nearly K\"ahler
structures based on the canonical almost complex structure on
homogeneous $3$-symmetric spaces (see \cite{S2}, \cite{WG},
\cite{G2}, \cite{Ki1}). It should be mentioned that the canonical
almost complex structure on such spaces became an effective tool
and a remarkable example in some deep constructions of
differential geometry and global analysis such as homogeneous
structures (\cite{TV}, \cite{Sat}, \cite{Ki4}, \cite{GV},
\cite{LV}, \cite{AG} etc.), Einstein metrics (\cite{SW},
\cite{SY}), holomorphic and minimal submanifolds (\cite{Sal1},
\cite{Sal2}), real Killing spinors (\cite{Gru}, \cite{BFGK},
\cite{Ka}).

The concept of generalized Hermitian geometry created in the 1980s
(see, for example, \cite{Ki2}, \cite{Ki7}) is a natural
consequence of the development of Hermitian geometry and the
theory of almost contact structures with many applications. One of
its central objects is the metric $f$-structures of the classical
type $(f^3+f=0)$, which include the class of almost Hermitian
structures. Many important classes of metric $f$-structures such
as K\"ahler, Killing, nearly K\"ahler, Hermitian $f$-structures
and some others were introduced and intensively investigated in
various aspects (see \cite{Ki2}, \cite{Ki5}, \cite{Ki7}, \cite{KL}
etc.). Specifically, Killing and nearly K\"ahler $f$-structures
became natural generalizations of classical nearly K\"ahler
structures in Hermitian geometry. However, this theory had not
provided new invariant examples of its own up to the recent
period, and so the lack of these examples was becoming all the
more noticeable.

There has recently been a qualitative change in the situation,
related to the complete solution of the problem of describing
canonical structures of classical type on regular $\Phi$-spaces
\cite{BS2}. A rich collection of canonical $f$-structures has been
discovered (including almost complex structures) leading to the
presentation of wide classes of invariant examples in generalized
Hermitian geometry (see \cite{B4}-\cite{B7}, \cite{C3} and
others). In particular, nearly K\"ahler $f$-structures were
provided with a remarkable class of their own invariant examples
(see \cite{B6}, \cite{B7}). This has ensured a continuation of the
classical results of J.A.Wolf, A.Gray, V.F.Kirichenko and others.
As to Killing $f$-structures, it is really an essential problem to
find proper non-trivial invariant examples of these structures.
Moreover, the possibilities for constructing such examples are
fairly limited (see \cite{B4}).

The main goals of this paper are

(i) to give a brief survey on invariant structures in generalized
Hermitian geometry and

(ii) to characterize all invariant $f$-structures on the flag
manifold $SU(3)/T_{max}$ in the sense of generalized Hermitian
geometry, in particular, to present first invariant examples of
Killing $f$-structures.

Sections 2-4 are mostly of survey character. In Section 2, we
collect some basic notions and results on homogeneous regular
$\Phi$-spaces and canonical affinor structures. In particular, a
precise description of all canonical structures of classical types
on homogeneous $k$-symmetric spaces is included. Besides, the
exact formulae for these structures and the relationship between
them on 4- and 5-symmetric spaces are presented.

In Section 3, we recall the main classes of almost Hermitian
structures following the Gray-Hervella division of almost
Hermitian manifolds into sixteen classes (see \cite{GH}). Besides,
we select particular results related to invariant almost Hermitian
structures.

Further, in Section 4, we describe main classes of metric
$f$-structures in generalized Hermitian geometry. Here we also
formulate the recent results on invariant nearly K\"ahler,
$G_1f$-, Hermitian, and Killing $f$-structures. In this
consideration, the canonical $f$-structures on homogeneous 4- and
5-symmetric spaces are especially important.

Finally, in Section 5, we examine in detail all invariant
$f$-structures on the complex flag manifold $SU(3)/T_{max}$ with
respect to all invariant Riemannian metrics. We discuss belonging
these structures to the main classes of metric $f$-structures
above mentioned. In particular, invariant non-trivial Killing
$f$-structures together with the corresponding Riemannian metrics
are first presented.

\section{\bf Homogeneous regular $\Phi$-spaces and canonical affinor structures}

Here we briefly formulate some basic definitions and results
related to regular $\Phi$-spaces and canonical affinor structures
on them. More detailed information can be found in \cite{BS2},
\cite{B10}, \cite{WG}, \cite{Ko}, \cite{F}, \cite{S1}, \cite{S2}.

Let $G$ be a connected Lie group, $\Phi$ its (analytic)
automorphism. Denote by $G^{\Phi}$ the subgroup of all fixed
points of $\Phi$ and $G_o^{\Phi}$ the identity component of
$G^{\Phi}$. Suppose a closed subgroup $H$ of $G$ satisfies the
condition $$G_o^{\Phi}\subset{H}\subset{G^{\Phi}}.$$ Then $G/H$ is
called  a {\it homogeneous $\Phi$-space}.

Homogeneous $\Phi$-spaces include homogeneous symmetric spaces
$(\Phi^2=id)$ and, more general, {\it homogeneous $\Phi$-spaces of
order $k$} $(\Phi^k=id)$ or, in the other terminology, {\it
homogeneous $k$-symmetric spaces} (see \cite{Ko}).

For any homogeneous $\Phi$-space $G/H$ one can define the mapping
$$
S_o = D\colon\ G/H \to G/H,\ xH\to \Phi (x) H.
$$
It is known \cite{S1} that $S_o$ is an analytic diffeomorphism of
$G/H$. $S_o$ is usually called a "symmetry" of $G/H$ at the point
$o=H$. It is evident that in view of homogeneity the "symmetry"
$S_p$ can be defined at any point $p\in G/H$. More exactly, for
any $p=\tau(x)o=xH,\ q=\tau(y)o=yH$ we put
$$
S_p=\tau(x)\circ S_o \circ \tau(x^{-1}).
$$
It is easy to show that
$$
S_p(yH)=x\Phi(x^{-1})\Phi(y)H.
$$
Thus any homogeneous $\Phi$-space is equipped with the set of
symmetries $\{S_p\mid p\in G/H\}$. Moreover, each $S_p$ is an
analytic diffeomorphism of the manifold $G/H$ (see \cite{S1}).

Note that there exist homogeneous $\Phi$-spaces that are not
reductive. That is why so-called regular $\Phi$-spaces first
introduced by N.A.Stepanov \cite{S1} are of fundamental
importance.

Let $G/H$ be a homogeneous $\Phi$-space, $\mathfrak{g}$ and
$\mathfrak{h}$ the corresponding Lie algebras for $G$ and $H$,
$\varphi=d{\Phi}_e$ the automorphism of $\mathfrak{g}$. Consider
the linear operator $A=\varphi-id$ and the Fitting decomposition
$\mathfrak{g}=\mathfrak{g}_0\oplus\mathfrak{g}_1$ with respect to
$A$, where   $\mathfrak{g}_0$ and  $\mathfrak{g}_1$ denote $0$-
and $1$-component of the decomposition respectively. Further, let
$\varphi=\varphi_{s}\:\varphi_{u}$ be the Jordan decomposition,
where $\varphi_{s}$ and $\varphi_{u}$ is a semisimple and
unipotent component of $\varphi$ respectively,
$\varphi_{s}\:\varphi_{u}=\varphi_{u}\:\varphi_{s}$. Denote by
$\mathfrak{g}^{\gamma}$ a subspace of all fixed points for a
linear endomorphism $\gamma$ in $\mathfrak{g}$. It is clear that
$\mathfrak{h}=\mathfrak{g}^{\varphi}=Ker\,A$,
$\mathfrak{h}\subset\mathfrak{g}_0$,
$\mathfrak{h}\subset\mathfrak{g}^{\varphi_s}$.

{\bf\noindent Definition 1} {\rm (\cite{S1}, \cite{BS2},
\cite{B10}, \cite{F}). A homogeneous $\Phi$-space $G/H$ is called
a {\it regular $\Phi$-space} if one of the following equivalent
conditions is satisfied:
\begin{enumerate}
\item $\mathfrak{h}=\mathfrak{g}_0$; \item
$\mathfrak{g}=\mathfrak{h}\oplus{A}\mathfrak{g}$; \item The
restriction of the operator $A$ to ${A}\mathfrak{g}$ is
non-singular; \item $A^2X=0\Longrightarrow{A}X=0$ for all
$X\in\mathfrak{g}$. \item The matrix of the automorphism $\varphi$
can be represented in the form \\ $\left(\begin{array}{cc}
  E & 0 \\
  0 & B
\end{array}\right),$ where the matrix $B$ does not admit the eigenvalue $1$.
\item $\mathfrak{h}=\mathfrak{g}^{\varphi_s}$.
\end{enumerate}}

\noindent We recall two basic facts:

\newpage

\begin{theorem}{\rm(\cite{S1})}
\begin{itemize}
\item Any homogeneous  $\Phi$-space of order $k$ $(\Phi^k\ = \
id)$ is a regular $\Phi$-space. \item Any regular $\Phi$-space is
reductive. More exactly, the Fitting decomposition

\begin{equation}\label{f1}
  \mathfrak{g}=\mathfrak{h}\oplus\mathfrak{m}, \
\mathfrak{m}=A\mathfrak{g}
\end{equation}
 is a reductive one.
\end{itemize}
\end{theorem}

Decomposition (\ref{f1})  is the {\it canonical reductive
decomposition} corresponding to a regular $\Phi$-space $G/H$, and
$\frak{m}$ is the {\it canonical reductive complement}.

We also note that for any regular $\Phi$-space $G/H$ each point
$p=xH\in G/H$ is an isolated fixed point of the symmetry $S_p$
(see \cite{S1}).

Decomposition (\ref{f1}) is obviously $\varphi$-invariant. Denote
by $\theta$ the restriction of $\varphi$ to $\mathfrak{m}$. As
usual, we identify $\mathfrak{m}$ with the tangent space
$T_o(G/H)$ at the point $o=H$. It is important to note that the
operator $\theta$ commutes with any element of the linear isotropy
group $Ad(H)$ (see \cite{S1}). It also should be noted (see
\cite{S1}) that
$$
(dS_o)_o=\theta.
$$

An {\it affinor structure} on a manifold is known to be a tensor
field of type $(1,1)$ or, equivalently, a field of endomorphisms
acting on its tangent bundle. Suppose $F$ is an invariant affinor
structure on a homogeneous manifold $G/H$. Then $F$ is completely
determined by its value $F_o$ at the point $o$, where $F_o$ is
invariant with respect to $Ad(H)$. For simplicity, we will denote
by the same manner both any invariant structure on $G/H$ and its
value at $o$ throughout the rest of the paper.

{\bf\noindent Definition 2} {\rm (\cite{BS1},\cite{BS2}). An
invariant affinor structure $F$ on a regular $\Phi$-space $G/H$ is
called {\em canonical} if its value at the point $o=H$ is a
polynomial in $\theta$.

Denote by $\mathcal A(\theta)$ the set of all canonical affinor
structures on a regular $\Phi$-space $G/H$. It is easy to see that
$\mathcal A(\theta)$ is a commutative subalgebra of the algebra
$\mathcal A$ of all invariant affinor structures on $G/H$.
Moreover, $$dim\ \mathcal A(\theta)=deg\ \nu\ \le\ dim\ G/H,$$
where $\nu$ is a minimal polynomial of the operator $\theta$. It
is evident that the algebra  $\mathcal A(\theta)$ for any
symmetric $\Phi$-space $(\Phi^2=id)$  consists of scalar
structures only, i.e. it is isomorphic to $\mathbb{R}$. As to
arbitrary regular $\Phi$-space $(G/H,\Phi)$, the algebraic
structure of its commutative algebra $\mathcal A(\theta)$ has been
recently completely described (see \cite{B9}).

It should be mentioned that all canonical structures are, in
addition, invariant with respect to the "symmetries" $\{S_p\}$ of
$G/H$(see \cite{S1}). Moreover, from $(dS_o)_o=\theta$ it follows
that the invariant affinor structure $F_p=(dS_p)_p, p\in G/H$
generated by the symmetries $\{S_p\}$ belongs to the algebra
$\mathcal A(\theta)$.

The most remarkable example of canonical structures is the
canonical almost complex structure
$J=\frac{1}{\sqrt{3}}(\theta-\theta^2)$ on a homogeneous
$3$-symmetric space (see \cite{S2}, \cite{WG}, \cite{G2}). It
turns out that it is not an exception. In other words, the algebra
$\mathcal A(\theta)$ contains many affinor structures of classical
types.

In the sequel we will concentrate on the following affinor
structures of classical types:

{\it almost complex structures} $J$ $(J^2=-1)$;

{\it almost product structures} $P$ $(P^2=1)$;

{\it $f$-structures} $(f^3+f=0)$ \cite{Y};

$f$-structures of hyperbolic type or, briefly, {\it
$h$-structures} $(h^3-h=0)$ \cite{Ki2}. \\ Clearly, $f$-structures
and $h$-structures are generalizations of structures $J$ and $P$
respectively.

All the canonical structures of classical type on regular
$\Phi$-spaces were completely described
\cite{BS1},\cite{BS2},\cite{B8}. In particular, for homogeneous
$k$-symmetric spaces, precise computational formulae were
indicated. For future reference we select here some results.

Denote by $\tilde s$ (respectively, $s$) the number of all
irreducible factors (respectively, all irreducible quadratic
factors) over $\mathbb{R}$ of a minimal polynomial $\nu$.

\begin{theorem}(\cite{BS1},\cite{BS2},\cite{B8})
Let $G/H$ be a regular $\Phi$-space.
\begin{enumerate}
  \item The algebra $\mathcal{A}(\theta )$ contains precisely $2^{\tilde s}$
structures $P$.
  \item $G/H$ admits a canonical structure $J$ if and only if
$s=\tilde{ s}$. In this case $\mathcal{ A}(\theta )$ contains
$2^s$ different structures $J$.
  \item $G/H$ admits a canonical $f$-structure if and only if
$s\ne 0$. In this case ${\mathcal A} (\theta )$ contains $3^s-1$
different $f$-structures. Suppose $s={\tilde s}$. Then $2^s$
$f$-structures are almost complex and the remaining $3^s-2^s-1$
have non-trivial kernels.
  \item The algebra $\mathcal{A}(\theta)$ contains $3^{\tilde{s}}$
  different $h$-structures. All these structures form a
  (commutative) semigroup in $\mathcal{A}(\theta)$ and include a
  subgroup of order $2^{\tilde{s}}$ of canonical structures $P$.
\end{enumerate}
\end{theorem}

Further, let $G/H$ be a homogeneous $k$-symmetric space. Then
$\tilde{s}=s+1$ if $-1\in \, spec\,\theta$, and $\tilde{s}=s$ in
the opposite case. We indicate explicit formulae enabling us to
compute all canonical $f$-structures and $h$-structures. We shall
also use the notation $$u= \left\{\begin{array}{ccc}
  n & {\text if} & k=2n+1 \\
  n-1 & {\text if} & k=2n
\end{array}\right..$$
\begin{theorem}(\cite{BS1},\cite{BS2},\cite{B8})
Let $G/H$ be a homogeneous $\Phi$-space of order $k$.
\begin{enumerate}
\item All non-trivial canonical $f$-structures on $G/H$ can be
given by the operators
$$f=\frac{2}{k}\sum_{m=1}^u\left(\sum_{j=1}^u\zeta_j\sin{
\frac{2\pi m j}{k}}\right)\left(\theta^m-\theta^{k-m}\right),$$
where $\zeta_j\in\{-1;0 ;1\},\;j=1,2,\ldots,u$, and not all
coefficients $\zeta_j$ are zero. In particular, suppose that
$-1\notin \, spec\,\theta$. Then the polynomials $f$ define
canonical almost complex structures $J$ iff all
$\zeta_j\in\{-1;1\}$. \item All canonical $h$-structures on $G/H$
can be given by the polynomials
$h=\sum\limits_{m=0}^{k-1}a_m\theta^m$, where:
\begin{enumerate}
\item if $k=2n+1$, then $$a_m=a_{k-m}=\frac{2}{k}\sum_{j=1}^u\xi_j
\cos{\frac{2\pi m j}{k}};$$ \item if $k=2n$, then $$a_m=a_{k-m}=
\frac{1}{k}\left(2\sum_{j=1}^u\xi_j \cos{\frac{2\pi m j}{k}} +
(-1)^m\xi_n \right)$$
\end{enumerate}

Here the numbers $\xi_j$ take their values from the set $\{-1;0
;1\}$ and the polynomials $h$ define canonical structures $P$ iff
all $\xi_j\in\{-1;1\}$.
\end{enumerate}
\end{theorem}

We now particularize the results above mentioned for homogeneous
$\Phi$-spaces of orders $3$, $4$, and $5$ only. Note that there
are no fundamental obstructions to considering of higher orders
$k$.
\begin{corol}\label{c1}(\cite{BS2},\cite{B8})
Let $G/H$ be a homogeneous $\Phi$-space of order $3$. There are
(up to sign) only the following canonical structures of classical
type on $G/H$: $$J=\frac{1}{\sqrt{3}}(\theta-\theta^2),\;P=1.$$
\end{corol}
We note that the existence of the structure $J$ and its properties
are well known (see \cite{S2},\cite{WG},\cite{G2},\cite{Ki1}).
\begin{corol}\label{c2}(\cite{BS2},\cite{B8})
On a homogeneous $\Phi$-space of order $4$ there are (up to sign)
the following canonical classical structures:
$$P=\theta^2,\;f=\frac12(\theta-
\theta^3),\;h_1=\frac12(1-\theta^2),\;h_2=\frac12(1+\theta^2).$$
The operators $h_1$ and $h_2$ form a pair of complementary
projectors: $h_1+h_2=1$, $h_1^2=h_1$, $h_2^2=h_2$. Moreover, the
following conditions are equivalent:
\begin{enumerate}
\item $-1\notin spec\,\theta$; \item the structure $P$ is trivial
$P=-1$; \item the $f$-structure is an almost complex structure;
\item the structure $h_1$ is trivial $(h_1=1)$; \item the
structure $h_2$ is null.
\end{enumerate}
\end{corol}

General properties of the canonical structures $P$ and $f$ on
homogeneous $4$-symmetric spaces were investigated in \cite{BD}.

\begin{corol}\label{c3}(\cite{BS1},\cite{BS2},\cite{B8})
There exist (up to sign) only the following   canonical structures
of classical type on any  homogeneous $\Phi$-space of order $5$:
\begin{gather*}
P=\frac{1}{\sqrt{5}}(\theta-\theta^2-\theta^3+\theta^4);\\
J_1=\alpha(\theta-\theta^4)-\beta(\theta^2-\theta^3);\quad
J_2=\beta(\theta-\theta^4)+\alpha(\theta^2-\theta^3);\\
f_1=\gamma(\theta-\theta^4)+\delta(\theta^2-\theta^3);\quad
f_2=\delta(\theta-\theta^4)-\gamma(\theta^2-\theta^3);\\
h_1=\frac12(1+P);\quad h_2=\frac12(1-P);
\end{gather*}
where $\alpha=\frac{\sqrt{5+2\sqrt5}}{5}$;\
$\beta=\frac{\sqrt{5-2\sqrt5}}{5}$;\
$\gamma=\frac{\sqrt{10+2\sqrt5}}{10}$;\
$\delta=\frac{\sqrt{10-2\sqrt5}}{10}$.

Besides, the following relations are satisfied:
\begin{gather*}
J_1\,P=J_2;\quad f_1\,P=J_1\,h_1=J_2\,h_1=f_1;\quad h_1\,P=h_1;\quad h_2\,P=-h_2;\\
 f_2\,P=J_2\,h_2=-J_1\,h_2=-f_2;\quad f_1\,f_2=h_1\,h_2=0;\quad h_1+h_2=P.
\end{gather*}

In addition, the following conditions are equivalent:
\begin{enumerate}
\item $spec\ \theta$ consists of two elements; \item the structure
$P$ is trivial; \item the structures $J_1$ and  $J_2$ coincide (up
to sign); \item one of the structures $f_1$ and $f_2$ is null,
while the other is an almost complex structure coinciding with one
of the structures $J_1$ and  $J_2$; \item one of the structures
$h_1$ and $h_2$ is trivial, while the other is null.
\end{enumerate}
\end{corol}

We note that for the first time the canonical structure $P$ on
homogeneous $5$-symmetric spaces was introduced and studied in
\cite{BC}. Other canonical structures on these spaces were later
studied in \cite{C1}-\cite{C3}.

It should be also mentioned that in the particular case of
homogeneous $\Phi$-spaces of any odd order $k=2n+1$ the method of
constructing invariant almost complex structures was described in
\cite{Ko}. It can be easily seen that all these structures are
canonical in the above sense.

\section{\bf Almost Hermitian structures}

We briefly recall some notions of Hermitian geometry including the
main classes of almost Hermitian structures.

Let $M$ be a smooth manifold, $\frak{X} (M)$ the Lie algebra of
all smooth vector fields on $M$, $d$ the exterior differentiation
operator. An {\it almost Hermitian structure} on $M$ (briefly,
{\it $AH$-structure}) is a pair $(g,J)$, where
$g=\langle\cdot,\cdot\rangle$ is a (pseudo)Riemannian metric on
$M$, $J$ an almost complex structure such that $\langle JX, JY
\rangle = \langle X, Y \rangle$ for any $X,Y\in\frak{X}(M)$. It
follows immediately that the tensor field $\Omega (X,Y)=\langle X,
JY \rangle$ is skew-symmetric, i.e. $(M,\Omega)$ is an almost
symplectic manifold. $\Omega$ is usually called a {\it fundamental
form} (the {\it K\"{a}hler form}) of an $AH$-structure $(g,J)$.

Further, denote by $\nabla$ the Levi-Civita connection of the
metric $g$ on $M$. We recall below some main classes of
$AH$-structures together with their defining properties (see, for
example, \cite {GH}):

\begin{tabbing}
  Kill f123 \= $f$-structure of class $G_1$, or -structure \= $T(X,X)=0,$ i.e.
  $\frak{X}(M)$ is an anticommutative $Q$-algebra \kill
  {\bf K} \> {\it K\"{a}hler structure}: \> $\nabla J=0$; \\
  {\bf H} \> {\it Hermitian structure}: \> $\nabla_X(J)Y-\nabla_{JX}(J)JY=0$; \\
  {\bf G$_1$} \> {$AH$-structure of class $G_1$,} or  \>
  $\nabla_X(J)X-\nabla_{JX}(J)JX=0$; \\
       \> $G_1$-{\it structure}: \> \\
  {\bf QK} \> {\it quasi-K\"{a}hler structure}: \> $\nabla_X(J)Y+\nabla_{JX}(J)JY=0$;\\
  {\bf AK} \> {\it almost K\"{a}hler structure}: \> $d\,\Omega=0$;\\
  {\bf NK} \> {\it nearly K\"{a}hler structure,}
  \> $\nabla_{X}(J)X=0$.\\
     \> or $NK$-{\it structure}: \>
\end{tabbing}

It is well known (see, for example, \cite {GH}) that
$$
{\bf K}\subset{\bf H}\subset{\bf G_1};\;\;{\bf K}\subset{\bf
NK}\subset{\bf G_1};\;\;{\bf NK}={\bf G_1}\cap{\bf QK};\;\;{\bf
K}={\bf H}\cap{\bf QK}.
$$

As usual, we will denote by $N$ the Nijenhuis tensor of an almost
complex structure $J$, that is,
$$
N(X,Y)=\frac14([JX,JY]-J[JX,Y]-J[X,JY]-[X,Y])
$$
for any $X,Y\in\frak{X}(M)$. Then the condition $N=0$ is
equivalent to the integrability of $J$ (see, for example, \cite
{KN}). Moreover, an almost Hermitian structure $(g,J)$ belongs to
the class {\bf H} if and only if $N=0$ (see, for example, \cite
{GH}).

As was already mentioned, the role of homogeneous almost Hermitian
manifolds is particularly important "because they are the model
spaces to which all other almost Hermitian manifolds can be
compared" (see \cite {G3}). A wealth of examples for the most
classes above noted, both of general and specific character, can
be found in \cite {WG}, \cite {G2}, \cite {G3}, \cite {Ki1} and
others. In particular, after the detailed investigation of the
6-dimensional homogeneous nearly K\"{a}hler manifolds
V.F.Kirichenko proved \cite {Ki1} that naturally reductive
strictly nearly K\"{a}hler manifolds $SO(5)/U(2)$ and
$SU(3)/T_{max}$ are not isometric even locally to the
6-dimensional sphere $S^6$. These examples gave a negative answer
to the conjecture of S.Sawaki and Y.Yamanoue (see \cite {SaYa})
claimed that any $6$-dimensional strictly $NK$-manifold was a
space of constant curvature. It should be noted that the canonical
almost complex structure $J=\frac1{\sqrt{3}}(\theta-\theta^2)$ on
homogeneous $3$-symmetric spaces plays a key role in these and
other examples of homogeneous $AH$-manifolds.

Let $G$ be a connected Lie group, $H$ its closed subgroup, $g$ an
invariant (pseudo-)Riemannian metric on the homogeneous space
$G/H$. Denote by $\mathfrak{g}$ and  $\mathfrak{h}$ the Lie
algebras corresponding to $G$ and $H$ respectively. Suppose that
$G/H$ is a reductive homogeneous space,
$\mathfrak{g}=\mathfrak{h}\oplus\mathfrak{m}$ the reductive
decomposition of the Lie algebra $\mathfrak{g}$. As usual, we
identify $\mathfrak{m}$ with the tangent space $T_o(G/H)$ at the
point $o=H$. Then the invariant metric $g$ is completely defined
by its value at the point $o$. For convenience we denote by the
same manner both any invariant metric $g$ on $G/H$ and its value
at $o$.

Recall that $(G/H,g)$ is {\it naturally reductive} with respect to
a reductive decomposition
$\mathfrak{g}=\mathfrak{h}\oplus\mathfrak{m}$  \cite{KN} if
$$g([X,Y]_{\mathfrak{m}},Z)=g(X,[Y,Z]_{\mathfrak{m}})$$ for all
$X,Y,Z\in\mathfrak{m}$. Here the subscript $\mathfrak{m}$ denotes
the projection of $\mathfrak{g}$ onto $\mathfrak{m}$ with respect
to the reductive decomposition.

We select here some known results closely related to the main
subject of our future consideration.

\begin{theorem} \label{t1}(\cite{AG})
Any invariant almost Hermitian structure on a naturally reductive
space $(G/H,g)$ belongs to the class {\bf G$_1$}.
\end{theorem}

\begin{theorem} \label{t2}(\cite{WG}, \cite{G2})
A homogeneous $3$-symmetric space $G/H$ with the canonical almost
complex structure $J$ and an invariant compatible metric $g$ is a
quasi-K\"{a}hler manifold. Moreover, $(G/H,J,g)$ belongs to the
class {\bf NK} if and only if $g$ is naturally reductive.
\end{theorem}

\begin{theorem} \label{t3}(\cite{Ma}, \cite{G4}, \cite{Ki6})
A $6$-dimensional strictly nearly K\"{a}hler manifold is Einstein.
\end{theorem}

Note that the latter result was obtained in \cite{Ki6} as a
particular case of the general approach based on investigating
nearly K\"{a}hler manifolds of constant type.

\section{\bf Metric $f$-structures and homogeneous manifolds}

The concept of the {\it generalized Hermitian geometry} created in
the 1980s (see, for example, \cite{Ki2}, \cite{Ki7}, \cite{Ki8})
was a natural consequence of the development of Hermitian geometry
and the theory of almost contact metric structures with many
applications. One of its central objects is the metric
$f$-structures of classical type, which include the class of
almost Hermitian structures. We recall here some basic notions.

An {\it $f$-structure} on a manifold $M$ is known to be a field of
endomorphisms $f$ acting on its tangent bundle and satisfying the
condition $f^3+f=0$ (see \cite{Y}). The number $r=dim\;Im\:f$ is
constant at any point of $M$ (see \cite{St}) and called a {\it
rank} of the $f$-structure. Besides, the number
$dim\;Ker\:f=dim\:M-r$ is usually said to be a {\it deficiency} of
the $f$-structure and denoted by $def\:f$.

Recall that an $f$-structure on a (pseudo)Riemannian manifold
$(M,g=\langle\cdot,\cdot\rangle)$ is called a {\it metric
$f$-structure}, if $\langle {fX}, Y \rangle +\langle X, fY \rangle
=0$, \;$X,Y\in\frak{X}(M)$ (see \cite{Ki2}). In the case the
triple $(M,g,f)$ is called a {\it metric $f$-manifold}. It is
clear that the tensor field $\Omega(X,Y)=\langle X, fY \rangle$ is
skew-symmetric, i.e. $\Omega$ is a $2$-form on $M$. $\Omega$ is
called a {\it fundamental form} of a metric $f$-structure
\cite{Ki7}, \cite{Ki2}. It is easy to see that the particular
cases $def\;f=0$ and $def\; f=1$ of metric $f$-structures lead to
almost Hermitian structures and almost contact metric structures
respectively.

Let $M$ be a metric $f$-manifold. Then
$\frak{X}(M)=\mathcal{L}\oplus\mathcal{M}$, where
$\mathcal{L}=Im\;f$ and $\mathcal{M}=Ker\;f$ are mutually
orthogonal distributions, which are usually called the {\it first}
and the {\it second fundamental distributions}  of the
$f$-structure respectively. Obviously, the endomorphisms $l=-f^2$
and $m=id+f^2$ are mutually complementary projections on the
distributions $\mathcal{L}$ and $\mathcal{M}$ respectively. We
note that in the case when the restriction of $g$ to $\mathcal{L}$
is non-degenerate the restriction $(F,g)$ of a metric
$f$-structure to $\mathcal{L}$ is an almost Hermitian structure,
i.e. $F^2=-id,\;\langle FX, FY \rangle=\langle X, Y
\rangle,\;X,Y\in\mathcal{L}$.

A fundamental role in the geometry of metric $f$-manifolds is
played by the {\it composition tensor} $T$, which was explicitly
evaluated in \cite{Ki7}:
\begin{equation}\label{eqT}
T(X,Y)=\frac14{f}(\nabla_{fX}(f){fY}-\nabla_{f^2X}(f){f^2Y}),
\end{equation}
where $\nabla$ is the Levi-Civita connection of a
(pseudo)Riemannian manifold $(M,g)$, \; $X,Y\in\frak{X}(M)$.

Using this tensor $T$, the algebraic structure of a so-called {\it
adjoint $Q$-algebra} in $\frak{X}(M)$ can be defined by the
formula:
$$
X\ast{Y}=T(X,Y).
$$
It gives the opportunity to introduce some classes of metric
$f$-structures in terms of natural properties of the adjoint
$Q$-algebra (see \cite{Ki2}, \cite{Ki7}).

We enumerate below the main classes of metric $f$-structures
together with their defining properties:

\begin{tabbing}
  Kill f123 \= $f$-structure of class $G_1$, or -structure \= $T(X,X)=0,$ i.e.
  $\frak{X}(M)$ is an anticommutative $Q$-algebra \kill
  {\bf Kf} \> {\it K\"{a}hler $f$--structure}: \> $\nabla f=0$; \\
  {\bf Hf} \> {\it Hermitian $f$--structure}: \> $T(X,Y)=0,$ i.e.
  $\frak{X}(M)$ is\\
     \>    \> an abelian $Q$-algebra;\\
  {\bf G$_1$f} \> {$f$-structure of class $G_1$,} or  \> $T(X,X)=0,$
  i.e.
  $\frak{X}(M)$ is \\
       \> $G_1f$-{\it structure}: \> an anticommutative $Q$-algebra;\\
  {\bf QKf} \> {\it quasi-K\"{a}hler $f$--structure}: \> $\nabla_X f +T_X
  f=0$;\\
  {\bf Kill f} \> {\it Killing $f$-structure}: \> $\nabla_X (f) X=0$;\\
  {\bf NKf} \> {\it nearly K\"{a}hler $f$-structure,}
  \> $\nabla_{fX}(f)fX=0$.\\
     \> or $NKf$-{\it structure}: \>
\end{tabbing}

The classes {\bf Kf}, {\bf Hf}, {\bf G$_1$f}, {\bf QKf} (in more
general situation) were introduced in \cite{Ki2} (see also
\cite{Si}). Killing $f$-manifolds {\bf Kill f} were defined and
studied in \cite{Gr1}, \cite{Gr2}. The class {\bf NKf} was first
determined in \cite{B1} (see also \cite{B6}, \cite{B7}).

The following relationships between the classes mentioned are
evident:
$$
{\bf Kf}={\bf Hf}\cap{\bf QKf};\;\;\; {\bf Kf}\subset{\bf
Hf}\subset{\bf G_1f};\;\;\; {\bf Kf}\subset{\bf
Kill\;f}\subset{\bf NKf}\subset{\bf G_1f.}
$$
It is important to note that in the special case $f=J$ we obtain
the corresponding classes of almost Hermitian structures (see
\cite{GH}). In particular, for $f=J$ the classes {\bf Kill f} and
{\bf NKf} coincide with the well-known class {\bf NK} of {\it
nearly K\"ahler structures}.

{\bf Remark 1.} Killing $f$-manifolds are often defined by
requiring the fundamental form $\Omega$ to be a Killing form, i.e.
$d\Omega=\nabla\Omega$ (see \cite{Gr1}, \cite{KL}). It is not hard
to prove that the definition is equivalent to the above condition
$\nabla_{X}(f)X=0$.

Indeed, in accordance with \cite{YB}, $\Omega$ is a Killing
$2$-form if and only if $\nabla\Omega$ is a $3$-form. Further,
using \cite{KN}, we have
\begin{equation}\label{eq1}
(\nabla\Omega)(X,Y;Z)=Z\:\Omega(X,Y)-\Omega(\nabla_{Z}X,Y)-\Omega(X,\nabla_{Z}Y).
\end{equation}
It follows that $\nabla\Omega$ is always skew-symmetric in
arguments $X$ and $Y$. Using the fact, it is easy to see that
$\nabla\Omega$ is a $3$-form if and only if $\nabla\Omega$ is
skew-symmetric in $Y$ and $Z$:
\begin{equation}\label{eq2}
(\nabla\Omega)(X,Y;Z)=-(\nabla\Omega)(X,Z;Y).
\end{equation}
Taking into account formula (\ref{eq1}) and the definition of
$\Omega$, condition (\ref{eq2}) can be written in the form:
$$
Z\langle X, fY \rangle-\langle \nabla_{Z}X, fY \rangle-\langle X,
f\nabla_{Z}Y \rangle+Y\langle X, fZ \rangle-\langle \nabla_{Y}X,
fYZ \rangle-\langle X, f\nabla_{Y}Z \rangle=0.
$$
Since $\nabla$ is the Levi-Civita connection, in particular, we
have:
$$
Z\langle X, Y \rangle=\langle \nabla_{Z}X, Y \rangle+\langle X,
\nabla_{Z}Y \rangle.
$$
It follows that the previous formula can be written in the form:
$$
\langle X, \nabla_{Z}fY -f\nabla_{Z}Y \rangle+\langle X,
\nabla_{Y}fZ -f\nabla_{Y}Z \rangle=0.
$$
Using the formula $\nabla_{X}(f)Y=\nabla_{X}fY-f\nabla_{X}Y$, we
get:
$$
\langle X, \nabla_{Z}(f)Y + \nabla_{Y}(f)Z \rangle=0.
$$
It implies that $\nabla_{Z}(f)Y + \nabla_{Y}(f)Z=0$ for any
$Y,Z\in\frak{X}(M)$. This is obviously equivalent to the condition
$\nabla_{X}(f)X=0,\;X\in\frak{X}(M)$. {\hfill $\square$}

Now we dwell on invariant metric $f$-structures on homogeneous
spaces.

Any invariant metric $f$-structure on a reductive homogeneous
space $G/H$ determines the orthogonal decomposition
$\mathfrak{m}=\mathfrak{m}_1\oplus\mathfrak{m}_2$ such that
$\mathfrak{m}_1=Im\:f$, $\mathfrak{m}_2=Ker\:f$.

As it was already noted (see Section 3), the main classes of
almost Hermitian structures are provided with the remarkable set
of invariant examples. It turns out that there is also a wealth of
invariant examples for the basic classes of metric $f$-structures.
These invariant metric $f$-structures can be realized on
homogeneous $k$-symmetric spaces with canonical $f$-structures. We
select here some results in this direction. More detailed
information can be found in  \cite{B1}-\cite{B7}, \cite{C3},
\cite{Li}.

\begin{theorem} \label{t4}(\cite{B5})
Any invariant metric $f$-structure on a naturally reductive space
$(G/H,g)$ is a $G_{1}f$-structure.
\end{theorem}

As a special case $(Ker\:f=0)$, it follows Theorem \ref{t1}.

\begin{theorem} \label{t5}(\cite{B5})
Let $(G/H,g,f)$ be a naturally reductive space with an invariant
metric $f$-structure that satisfies the condition
$[\frak{m}_1,\frak{m}_1]\subset\frak{m}_2\oplus\frak{h}$. Then
$G/H$ is a Hermitian $f$-manifold.
\end{theorem}

We note that Theorem \ref{t5} is also valid for arbitrary
invariant (pseudo)Riemanni\-an metric $g$ compatible with an
invariant $f$-structure on a reductive homogeneous space $G/H$
(see \cite{BV}).

Theorems \ref{t4} and \ref{t5} can be effectively provided with a
large class of examples. In particular, for a semi-simple group
$G$, the invariant  (pseudo)Riemannian metric $g$ generated by the
Killing form on any regular $\Phi$-space $G/H$ is naturally
reductive with respect to the canonical reductive decomposition
$\frak{g}=\frak{h}\oplus\frak{m}$ (see \cite{S1}). Moreover, all
canonical structures $f$ and $J$ on homogeneous naturally
reductive $k$-symmetric spaces are compatible with such a metric,
i.e. $f$ is a metric $f$-structure, $J$ is an almost Hermitian
structure (see \cite{B1}, \cite{B10}). In what follows, we mean by
a naturally reductive decomposition the canonical reductive
decomposition for a regular $\Phi$-space $G/H$. To sum up, we have

\begin{theorem} \label{t6}(\cite{B5})
Let $(G/H,g)$ be a naturally reductive $k$-symmetric space. Any
canonical metric $f$-structure on $G/H$ is a $G_{1}f$-structure,
and any canonical almost Hermitian structure $J$ is a
$G_{1}$-structure.
\end{theorem}

\begin{theorem} \label{t7}(\cite{B6},\cite{B7})
Let $G/H$ be a regular $\Phi$-space, $g$ a naturally reductive
metric on $G/H$ with respect to the canonical reductive
decomposition $\mathfrak{g}=\mathfrak{h}\oplus\mathfrak{m}$, $f$ a
metric canonical $f$-structure on $G/H$. Suppose the $f$-structure
satisfies the condition $f^2=\pm\theta\:f$. Then $(G/H,g,f)$ is a
nearly K\"ahler $f$-manifold.
\end{theorem}

\begin{corol} (\cite{B6},\cite{B7})
Let $(G/H,g)$ be a  naturally reductive homogeneous $\Phi$-space
of order $k=4n,\:n\ge{1}$. If $\{i,-i\}\ \subset\ spec\ \theta$,
then there exists a non-trivial canonical $NKf$-structure on
$G/H$.
\end{corol}

We stress the particular role of homogeneous $4$- and
$5$-symmetric spaces.

\begin{theorem} \label{t8}(\cite{B3}-\cite{B7})
The canonical $f$-structure $f=\frac12(\theta-\theta^3)$ on any
naturally reductive $4$-symmetric space $(G/H,g)$  is both a
Hermitian $f$-structure and a nearly K\"ahler $f$-structure.
Moreover, the following conditions are equivalent:

1) $f$ is a K\"ahler $f$-structure; 2) $f$ is a Killing
$f$-structure; 3) $f$ is a quasi-K\"ahler $f$-structure; 4) $f$ is
an integrable $f$-structure; 5)
$[\frak{m}_1,\frak{m}_1]\subset\frak{h}$; 6)
$[\frak{m}_1,\frak{m}_2]=0$; 7) $G/H$ is a locally symmetric
space: $[\frak{m},\frak{m}]\subset\frak{h}$.
\end{theorem}

\begin{theorem} \label{t9}(\cite{B2}-\cite{B5}, \cite{B7}, \cite{C3})
Let $(G/H,g)$ be a naturally reductive $5$-symmetric space, $f_1$
and $f_2$, $J_1$ and $J_2$ the canonical structures on this space.
Then $f_1$ and $f_2$ belong to both classes {\bf Hf} and {\bf
NKf}. Moreover, the following conditions are equivalent:

1) $f_1$ is a K\"ahler $f$-structure; 2) $f_2$ is a K\"ahler
$f$-structure; 3) $f_1$ is a Killing $f$-structure; 4) $f_2$ is a
Killing $f$-structure; 5) $f_1$ is a quasi-K\"ahler $f$-structure;
6) $f_2$ is a quasi-K\"ahler $f$-structure; 7) $f_1$ is an
integrable $f$-structure; 8) $f_2$ is an integrable $f$-structure;
9) $J_1$ and $J_2$ are $NK$-structures; 10)
$[\frak{m}_1,\frak{m}_2]=0$ (here $\frak{m}_1=Im\:f_1=Ker\:f_2,
\frak{m}_2=Im\:f_2=Ker\:f_1$); 11) $G/H$ is a locally symmetric
space: $[\frak{m},\frak{m}]\subset\frak{h}$.
\end{theorem}

It should be mentioned that Riemannian homogeneous $4$-symmetric
spaces of classical compact Lie groups were classified and
geometrically described in \cite{J}. The similar problem for
homogeneous $5$-symmetric spaces was considered in \cite{TX}. By
Theorem \ref{t8} and Theorem \ref{t9}, it presents a collection of
homogeneous $f$-manifolds in the classes {\bf NKf} and {\bf Hf}.
Note that the canonical $f$-structures under consideration are
generally non-integrable.

Besides, there are invariant $NKf$-structures and $Hf$-structures
on homogeneous spaces $(G/H,g)$, where the metric $g$ is not
naturally reductive. The example of such a kind can be realized on
the $6$-dimensional Heisenberg group $(N,g)$. As to details
related to this group, we refer to  \cite{Ka1}, \cite{Ka2},
\cite{TV}.

\begin{theorem} \label{t10}(\cite{B5}-\cite{B7})
The 6-dimensional generalized Heisenberg group $(N,g)$  with
respect to the canonical $f$-structure
$f=\frac12(\theta-\theta^3)$ of a homogeneous $\Phi$-space of
order $4$ is both $Hf$- and $NKf$-manifold. This $f$-structure is
neither Killing nor integrable on $(N,g)$.
\end{theorem}

{\bf Remark 2.} Theorems \ref{t8} and \ref{t10}, in particular,
illustrate simultaneously the analogy and the difference between
the canonical almost complex structure $J$ on homogeneous
$3$-symmetric spaces $(G/H,g,J)$ and the canonical $f$-structure
on homogeneous $4$-symmetric spaces $(G/H,g,f)$ (see Theorem
\ref{t2}).

Let us also remark that the 6-dimensional generalized Heisenberg
group $(N,g)$ is an example of solvable type. In Section 5, we
present $NKf$-structures with non-naturally reductive metrics in
the case of semi-simple type.

Finally, we briefly discuss the existence problem for invariant
Killing $f$-structu\-res. By Theorems \ref{t8} and \ref{t9}, the
canonical $f$-structures on naturally reductive $4$- and
$5$-symmetric spaces are never strictly (i.e. non-K\"ahler)
Killing $f$-structures. Moreover, we recall the following general
result:

\begin{theorem} \label{t11}(\cite{B4})
Let $(G/H,g,f)$ be a naturally reductive Killing $f$-manifold.
Then the following relations hold:
$$
[\mathfrak{m}_1,\mathfrak{m}_1]\subset\mathfrak{m}_1\oplus\mathfrak{h},
\quad
[\mathfrak{m}_2,\mathfrak{m}_2]\subset\mathfrak{m}_2\oplus\mathfrak{h},
\quad [\mathfrak{m}_1,\mathfrak{m}_2]\subset\mathfrak{h}.
$$
In particular, both the fundamental distributions of the Killing
$f$-structure generate invariant totally geodesic foliations on
$G/H$.
\end{theorem}

By the results in \cite{Gr1} and Theorem \ref{t11}, it follows

\begin{corol}(\cite{B4})
There are no non-trivial (i.e. $def\:f>0$) invariant Killing
$f$-structures of the so-called fundamental type (see \cite{Gr1})
on naturally reductive homogeneous spaces $(G/H,g)$.
\end{corol}

This fact is a wide generalization of the similar result of
A.Gritsans obtained for Riemannian globally symmetric spaces.
Besides, it shows a substantial difference between invariant
Killing $f$-structures and invariant $NK$-structures.

In Section 5, we will indicate, in particular, first examples of
invariant Killing $f$-structures.

\section{\bf Invariant $f$-structures on the complex \\flag manifold $\bf M=SU(3)/T_{max}$}

In this Section, we will consider all invariant $f$-structures on
the flag manifold $M=SU(3)/T_{max}$. Note that invariant almost
complex structures (i.e. $f$-structures of maximal rank $6$) on
this space were investigated in \cite{G3}, \cite{AGI1},
\cite{AGI2} and many other papers.

The homogeneous manifold $SU(3)/T_{max}$ is known to be an
important example in many branches of differential geometry and
beyond. In particular, $M=SU(3)/T_{max}$ is a Riemannian
homogeneous $3$-symmetric space not homeomorphic with the
underlying manifold $M$ of any Riemannian symmetric space (see
\cite{LP1}). Further, $M$ is a homogeneous $k$-symmetric space for
any $k\ge3$. Moreover, $M$ is a naturally reductive Riemannian
homogeneous space that is {\it non-commutative} (see \cite{J2}).
It means that the algebra of invariant differential operators
$\mathcal{D}(SU(3)/T_{max})$ is not commutative (see \cite{H1}).
It follows that $M=SU(3)/T_{max}$ is not even a {\it weakly
symmetric space} (see, for example, \cite{Vi}).

Besides, $M$ is the twistor space for the projective space
$\mathbb{C}P^2$ (see, for example, \cite{Be}, Chapter 13). It was
a key point for constructing the first examples of $6$-dimensional
Riemannian manifolds admitting a real Killing spinor (see
\cite{BFGK}). More exactly, the flag manifold $M=SU(3)/T_{max}$
with the nearly K\"ahler structure $(g,J)$ just possesses a real
Killing spinor (see \cite{BFGK}, \cite{Gru}). Moreover, using the
duality procedure for this space $SU(3)/T_{max}$, one can
effectively construct pseudo-Riemannian homogeneous manifolds with
the real Killing spinors (see \cite{Ka}).

Let $\Phi=I(s)$ be an inner automorphism of the Lie group $SU(3)$
defined by the element $s=diag\:(\varepsilon,\overline
\varepsilon,1)$, where $\varepsilon$ is a primitive third root of
unity. Then the subgroup $H=G^{\Phi}$ of all fixed points of
$\Phi$ is of the form:
$$
G^{\Phi}=\{ diag(e^{i\beta_1}, e^{i\beta_2},e^{i\beta_3})|
\beta_1+\beta_2+\beta_3=0,\;\beta_j\in\mathbb{R}\}.
$$
Obviously, $G^{\Phi}$ is isomorphic to $T^2=T_{max}$ diagonally
imbedded into $SU(3)$. It means that the flag manifold
$M=SU(3)/T_{max}$ is a homogeneous $3$-symmetric space defined by
the automorphism $\Phi$.

Consider the canonical reductive decomposition
$\frak{g}=\frak{h}\oplus\frak{m}$ of the Lie algebra
$\frak{g}=\frak{s}\frak{u}(3)$ for the homogeneous $\Phi$-space
$M$. Using the notations in \cite{R1}, we obtain:
$$\frak{g}=\frak{su}(3)=\left\{\left(\begin{array}{ccc}
  \alpha_1 & a & \overline{c} \\
  -\overline{a} & \alpha_2 & b \\
  -c & -\overline{b} & \alpha_3
\end{array}\right) \left| \begin{array}{l}
  \alpha_1,\alpha_2, \alpha_3\in Im\,\mathbb{C}, \\
  a,b,c\in\mathbb{C}, \\
  \alpha_1+\alpha_2+\alpha_3=0
\end{array}\right.\right\}=$$
$$=E(\alpha_1,\alpha_2, \alpha_3)\oplus
D(a,b,c)=\frak{h}\oplus\frak{m}.$$

If we put $X=D(a,b,c), Y=D(a_1,b_1,c_1),
Z=E(\alpha_1,\alpha_2,\alpha_3)$, then the Lie brackets can be
briefly indicated (see \cite{R2}):

$$ [X,Y]=D(\overline{bc_1-b_1c}, \overline{ca_1-c_1a},
\overline{ab_1-a_1b}) -$$  $$
2E(Im(a\overline{a_1}+\overline{c}c_1),Im(\overline{a}a_1+b\overline{b_1}),
Im(c\overline{c_1}+\overline{b}b_1) ) ,$$

$$[Z,X]=D(\alpha_1 a-a\alpha_2, \alpha_2 b-b\alpha_3,
\alpha_3c-c\alpha_1).$$

Further, we put
$\frak{m}=\frak{m}_1\oplus\frak{m}_2\oplus\frak{m}_3,$ where

$$\frak{m}_1 = \{X\in\frak{su}(3)| X=D(a,0,0), a\in\mathbb{C}\},$$
$$\frak{m}_2  =\{X\in\frak{su}(3)| X=D(0,b,0), b\in\mathbb{C}\},$$
$$\frak{m}_3 = \{X\in\frak{su}(3)| X=D(0,0,c), c\in\mathbb{C}\}.$$

Using the Killing form of the Lie algebra $\frak{su}(3)$, we
define an invariant inner product on $\frak{m}$:
$$
g_o(X,Y)=\langle X, Y\rangle_o=-\frac12Re\:tr\:XY.
$$
Then (see \cite{R1})
$\frak{g}=\frak{h}\oplus\frak{m}_1\oplus\frak{m}_2\oplus\frak{m}_3$
is $\langle\cdot,\cdot\rangle_o$-orthogonal decomposition
satisfying the following relations:
$$
[\frak{h},\frak{m}_j]\subset\frak{m}_j,\;[\frak{m}_j,\frak{m}_j]\subset\frak{h},\;
[\frak{m}_j,\frak{m}_{j+1}]\subset\frak{m}_{j+2},
$$
where $j=1,2,3$ and the index $j$ should be reduced by modulo $3$.
Besides, the $H$-modules $\frak{m}_j$ are pairwise non-isomorphic.

Now we turn to invariant Riemannian metrics on $M$. Taking into
account the well-known one-to-one correspondence between
$G$-invariant Riemannian metrics on $G/H$ and $Ad(H)$-invariant
inner products on $\frak{m}$ (see \cite{KN}), we will make use of
the following fact:

\begin{lemma}(\cite{R1})
Any $SU(3)$-invariant Riemannian metric
$g=\langle\cdot,\cdot\rangle$ on the flag manifold
$M=SU(3)/T_{max}$ can be written in the form
$$
g=\langle\cdot,\cdot\rangle=
\lambda_1\langle\cdot,\cdot\rangle_{o{|\frak{m}_1\times\frak{m}_1}}
+
\lambda_2\langle\cdot,\cdot\rangle_{o{|\frak{m}_2\times\frak{m}_2}}
+
\lambda_3\langle\cdot,\cdot\rangle_{o{|\frak{m}_3\times\frak{m}_3}},
$$
where $\lambda_j>0,\;j=1,2,3.$
\end{lemma}

A triple $(\lambda_1,\lambda_2,\lambda_3)$ is called \cite{R1} a
{\it characteristic collection} of a Riemannian metric $g$ above
mentioned . Considering Riemannian metrics up to homothety, one
can assume that
$(\lambda_1,\lambda_2,\lambda_3)=(1,t,s),\;t>0,s>0.$ For
convenience we will denote this correspondence in the following
way: $g=(\lambda_1,\lambda_2,\lambda_3)$ or $g=(1,t,s).$

We also recall the following result:
\begin{theorem} \label{t12}(\cite{ZW},\cite{AN},\cite{R1})
There are exactly (up to homothety) the following invariant
Einstein metrics on the flag manifold $SU(3)/T_{max}$ $:$
$$(1,1,1), (1,2,1), (1,1,2), (2,1,1).$$
\end{theorem}

Let $\alpha$ be the Nomizu function (see \cite{N}) of the
Levi-Civita connection $\nabla$ for an invariant Riemannian metric
$g=\langle\cdot,\cdot\rangle$ on a reductive homogeneous space
$G/H$. Then
\begin{equation}\label{eq4}
\alpha(X,Y)=\frac12[X,Y]_{\frak{m}}+U(X,Y),\;\;\;X,Y\in\frak{m},
\end{equation}
where $U:\frak{m}\times\frak{m}\to\frak{m}$ is a symmetric
bilinear mapping determined by the formula (see\cite{KN}):
$$
2\langle U(X,Y),Z\rangle=\langle X,[Z,Y]_{\frak{m}}\rangle+
\langle [Z,X]_{\frak{m}},Y\rangle.
$$

For our case in these notations we have
\begin{lemma}(\cite{W},\cite{R2})
For the Levi-Civita connection of a Riemannian metric
$g=(\lambda_1,\lambda_2,\lambda_3)$ on the flag manifold
$SU(3)/T_{max}$ the following conditions are satisfied$:$
$$
U(X,Y)=0, \;\;if\;\; X,Y\in\frak{m}_j,\,j\in\{1,2,3\};
$$
$$ U(X,Y)= -
(2\lambda_j)^{-1}(\lambda_{j+1}-\lambda_{j+2})[X,Y], \;\;if\;\;
X\in\frak{m}_{j+1},Y\in\frak{m}_{j+2},
$$
where $j=1,2,3$ and the numbers $j$ are reduced by modulo $3$.
\end{lemma}

Let us now turn to invariant $f$-structures on $M=SU(3)/T_{max}$.
Keeping the above notations, any invariant $f$-structure on $M$
can be expressed by the mapping
\begin{equation}\label{eqf}
f:D(a,b,c)\rightarrow D(\zeta_1ia,\zeta_2ib,\zeta_3ic),
\end{equation}
where $\zeta_{j}\in\{1,0,-1\},\;j=1,2,3$, $i$ is the imaginary
unit. We will call the collection $(\zeta_1,\zeta_2,\zeta_3)$ a
{\it characteristic collection} of the invariant $f$-structure and
for convenience denote $f=(\zeta_1,\zeta_2,\zeta_3).$ Obviously,
all invariant $f$-structures on $M$ pairwise commute.

If we agree to consider  $f$-structures up to sign, then there are
the following invariant $f$-structures on $M=SU(3)/T_{max}$:

1) {\it invariant $f$-structures of rank $6$ (invariant almost
complex structures)}:
$$
J_1=(1,1,1),\;\;J_2=(1,-1,1),\;\;J_3=(1,1,-1),\;\;J_4=(1,-1,-1).
$$

2) {\it invariant $f$-structures of rank $4$}:
$$
\begin{array}{ccc}
f_1=(1,1,0),\;\;\;\;f_2=(1,0,1),\;\;\;\;f_3=(0,1,1),\\
f_4=(1,-1,0),\;\;\;\;f_5=(1,0,-1),\;\;\;\;f_6=(0,1,-1).
\end{array}
$$

3) {\it invariant $f$-structures of rank $2$}:
$$
f_7=(1,0,0),\;\;\;\;f_8=(0,1,0),\;\;\;\;f_9=(0,0,1).
$$

Our description of all invariant $f$-structures and all invariant
Riemannian metrics evidently implies that any invariant
$f$-structure $f=(\zeta_1,\zeta_2,\zeta_3)$ is a metric
$f$-structure with respect to any invariant Riemannian metric
$g=(\lambda_1,\lambda_2,\lambda_3)$. In particular,
$J_j,\;j=1,2,3,4$ are invariant almost Hermitian structures with
respect to all invariant Riemannian metrics
$g=(\lambda_1,\lambda_2,\lambda_3)$.

Now we are able to investigate all invariant $f$-structures in the
sense of generalized Hermitian geometry, i.e. the special classes
{\bf Kf}, {\bf NKf}, {\bf Kill f}, {\bf Hf}, {\bf G$_1$f}.

A key point of our consideration belongs to the expression
$\nabla_X(f)Y$. Using formula (\ref{eq4}), we get:
$$
\begin{array}{cc}
\nabla_X(f)Y=\nabla_X fY-f\nabla_X Y=\alpha(X,fY)-f
\alpha(X,Y)\\
=\frac12 ([X,fY]_{\frak{m}}-f[X,Y]_{\frak{m}})+U(X,fY)-fU(X,Y).
\end{array}
$$
As a result, we can obtain:
\begin{multline}\label{123}
  \nabla_X(f)Y=\frac12 D(A,B,C), \text{where}\\
  A=\overline{i((\zeta_1+\zeta_3)(1+s-t)bc_1+(\zeta_1+\zeta_2)(s-t-1)b_1c)},\\
  B=\overline{i((\zeta_2+\zeta_1)(1+\frac{1-s}{t})ca_1+(\zeta_2+\zeta_3)(\frac{1-s}{t}-1)c_1a)},\\
  C=\overline{i((\zeta_3+\zeta_2)(\frac{t-1}{s}+1)ab_1+(\zeta_3+\zeta_1)(\frac{t-1}{s}-1)a_1b)}.
\end{multline}

\subsection{\bf K\"ahler $f$-structures}

K\"ahler $f$-structures are defined by the condition
$\nabla_X(f)Y=0$ (see Section 4). Using formula (\ref{123}), this
condition is equivalent to the following system of equations:
\begin{equation}\label{syst1}
  \left\{\begin{array}{c}
    (\zeta_1+\zeta_3)(s-t+1)=0 \\
    (\zeta_1+\zeta_2)(s-t-1)=0 \\
    (\zeta_2+\zeta_3)(s+t-1)=0 \
  \end{array}\right.
\end{equation}
Solving system (\ref{syst1}) for all invariant $f$-structures, we
obtain the following result:

\begin{propos}
The flag manifold $M=SU(3)/T_{max}$ admits the following invariant
K\"ahler $f$-structures with respect to the corresponding
invariant Riemannian metrics only:
$$
\begin{array}{lll}
J_2=(1,-1,1),&g_t=(1,t,t-1),&t>1;\\
J_3=(1,1,-1),&g_t=(1,t,t+1),&t>0;\\
J_4=(1,-1,-1),&g_t=(1,t,1-t),&0<t<1.
\end{array}
$$
In particular, there are no invariant K\"ahler $f$-structures of
rank $2$ and $4$ on $M$.
\end{propos}

We note that the result is known for invariant almost complex
structures (see \cite{G3},\cite{AGI2}). We can also observe that
for each of K\"ahler $f$-structures $J_2,J_3,J_4$ the
corresponding $1$-parameter set $g_t$ of invariant Riemannian
metrics contains exactly one Einstein metric excluding the
naturally reductive metric $g=(1,1,1)$ (see Theorem \ref{t12}).
Taking into account Theorem \ref{t2}, the latter fact implies that
the structures $J_2,J_3,J_4$ cannot be realized as the canonical
almost complex structures $J=\frac1{\sqrt3}(\theta-\theta^2)$ for
some homogeneous $\Phi$-spaces of order $3$.

In addition, Lie brackets relations for the subspaces $\frak{m}_j,
\;j=1,2,3$ imply that all invariant $f$-structures of rank $2$ and
$4$ are non-integrable. It immediately follows that these
$f$-structures cannot be K\"ahler $f$-structures.

\subsection{\bf Killing $f$-structures}

The defining condition for Killing $f$-structures can be written
in the form $\nabla_X(f)X=0$ (see Section 4). From (\ref{123}), it
follows
\begin{multline*}
  \nabla_X(f)X=\frac12 D(A_0,B_0,C_0), \text{where}\\
  A_0=\overline{ibc((\zeta_1+\zeta_3)(1+s-t)+(\zeta_1+\zeta_2)(s-t-1))},\\
  B_0=\overline{ica((\zeta_2+\zeta_1)(1+\frac{1-s}{t})+(\zeta_2+\zeta_3)(\frac{1-s}{t}-1))},\\
  C_0=\overline{iab((\zeta_3+\zeta_2)(\frac{t-1}{s}+1)+(\zeta_3+\zeta_1)(\frac{t-1}{s}-1))}.
\end{multline*}
It easy to show that the condition $\nabla_X(f)X=0$ is equivalent
to the following system of equations:
$$
  \left\{\begin{array}{c}
    (\zeta_1+\zeta_3)(s-t+1)+(\zeta_1+\zeta_2)(s-t-1)=0 \\
    (\zeta_1+\zeta_2)(s-t-1)+(\zeta_2+\zeta_3)(s+t-1)=0 \
  \end{array}\right.
$$
Analyzing this system for all invariant $f$-structures, we obtain
the following result:

\begin{propos}
All invariant strictly Killing (i.e. non-K\"ahler) $f$-structures
on the flag manifold $M=SU(3)/T_{max}$ and the corresponding
invariant Riemannian metrics (up to homothety) are indicated
below:
$$
\begin{array}{ll}
J_1=(1,1,1),&g=(1,1,1);\\
f_1=(1,1,0),&g=(3,3,4);\\
f_2=(1,0,1),&g=(3,4,3);\\
f_3=(0,1,1),&g=(4,3,3).\
\end{array}
$$
In particular, there are no invariant Killing $f$-structures of
rank $2$ on $M$.
\end{propos}

Note the structure $J_1$ is a well-known non-integrable nearly
K\"ahler structure on a naturally reductive space $M$ (see
\cite{G2}, \cite{G3}, \cite{Ki1}, \cite{AGI2} and others). The
structures $f_1,f_2,f_3$ present first invariant non-trivial
Killing $f$-structures \cite{B11}. The important feature of these
structures is that the corresponding invariant Riemannian metrics
are not Einstein (see Theorem \ref{t12}). It illustrates a
substantial difference between non-trivial strictly Killing
$f$-structures and strictly $NK$-structures at least in the
$6$-dimensional case (see Theorem \ref{t3}).

{\bf Remark 3.} It is interesting to note that all strictly
Killing $f$-structures above indicated are canonical
$f$-structures for suitable homogeneous $\Phi$-spaces of the Lie
group $SU(3)$. We already mentioned that $M=SU(3)/T_{max}$ is a
homogeneous $k$-symmetric space for any $k\ge 3$. It means $M$ as
an underlying manifold could be generated by various automorphisms
$\Phi$ of the Lie group $SU(3)$. In particular, $J_1$ is the
canonical almost complex structure
$J=\frac1{\sqrt{3}}(\theta-\theta^2)$ for the homogeneous
$\Phi$-space of order $3$, where $\Phi=I(s),\;
s=diag\:(\varepsilon,\overline
\varepsilon,1),\;\varepsilon=\sqrt[3]{1}$ (see the beginning of
this Section). Further, if we consider the automorphism
$\Phi_1=I(s_1),\; s_1=diag\:(i,-i,1)$, where $i=\sqrt[4]{1}$ is
the imaginary unit, then $M$ is a homogeneous $\Phi_1$-space of
order $4$. The corresponding canonical $f$-structure
$f=\frac12(\theta_1-\theta_{1}^3)$ for this $\Phi_1$-space just
coincides (up to sign) with the $f$-structure $f_3=(0,1,1)$. The
structures $f_1$ and $f_2$ can be obtained in the similar way.
Moreover, all the structures $f_1,f_2,f_3$ and $f_7,f_8,f_9$ can
be realized as canonical $f$-structures for suitable homogeneous
$\Phi$-spaces of order $5$.

We also note that all $f$-structures $f_1,f_2,f_3$ are just the
restrictions of the structure $J_1$ onto the corresponding
distributions $\frak{m}_p\oplus\frak{m}_q,\;p,q\in\{1,2,3\}.$

\subsection{\bf Nearly K\"ahler $f$-structures}

Using (\ref{123}), we can easily obtain:
\begin{multline*}
  \nabla_{fX}(f)fX=\frac12 D(\hat A,\hat B,\hat C), \text{where}\\
  \hat A=\overline{-i\zeta_2\zeta_3bc((\zeta_1+\zeta_3)(1+s-t)+(\zeta_1+\zeta_2)(s-t-1))},\\
  \hat B=\overline{-i\zeta_1\zeta_3ca((\zeta_2+\zeta_1)(1+\frac{1-s}{t})+
  (\zeta_2+\zeta_3)(\frac{1-s}{t}-1))},\\
  \hat C=\overline{-i\zeta_1\zeta_2ab((\zeta_3+\zeta_2)(\frac{t-1}{s}+1)+
  (\zeta_3+\zeta_1)(\frac{t-1}{s}-1))}.
\end{multline*}
It follows that the condition $\nabla_{fX}(f)fX=0$ is reduced to
the following system of equations:
$$
  \left\{\begin{array}{c}
    \zeta_2\zeta_3((\zeta_1+\zeta_3)(s-t+1)+(\zeta_1+\zeta_2)(s-t-1))=0 \\
    \zeta_1\zeta_3((\zeta_2+\zeta_1)(1+t-s)+(\zeta_2+\zeta_3)(1-s-t))=0 \\
    \zeta_1\zeta_2((\zeta_3+\zeta_2)(t+s-1)+(\zeta_3+\zeta_1)(t-s-1))=0 \
  \end{array}\right.
$$
Consideration of this system implies

\begin{propos}\label{p3}
The only invariant strictly nearly K\"ahler (i.e. non-K\"ahler)
$f$-structure of rank $6$ on the flag manifold $M=SU(3)/T_{max}$
is the nearly K\"ahler structure $J_1=(1,1,1)$ with respect to the
naturally reductive metric $g=(1,1,1)$.

Invariant strictly nearly K\"ahler $f$-structures of rank $4$ and
the corresponding invariant Riemannian metrics (up to homothety)
on $M$ are:
$$
\begin{array}{lll}
f_1=(1,1,0),&g_s=(1,1,s),&s>0;\\
f_2=(1,0,1),&g_t=(1,t,1),&t>0;\\
f_3=(0,1,1),&g_t=(1,t,t),&t>0.\
\end{array}
$$
The invariant $f$-structures $f_7,f_8,f_9$ of rank $2$ on $M$ are
strictly $NKf$-structures with respect to all invariant Riemannian
metrics $g=(1,t,s),\;t,s>0$.
\end{propos}

First we notice that the structures $f_1,f_2,f_3$ and
$f_7,f_8,f_9$ provide invariant examples of $NKf$-structures with
non-naturally reductive metrics on the homogeneous space
$M=SU(3)/T_{max}$, which belongs to a semi-simple type.

We can also observe that for any invariant strictly
$NKf$-structure on $M$ there exists at least one (up to homothety)
corresponding Einstein metric. More exactly, for these
$NKf$-structures of rank $6$, $4$, and $2$ there are (up to
homothety) $1$, $2$, and $4$ Einstein metrics respectively (see
Theorem \ref{t12}). In a certain degree, it is a particular
analogy with the result of Theorem \ref{t3}. This particular fact
and some related general results lead to the following conjecture,
which seems to be plausible:

{\bf Conjecture}. For any strictly nearly K\"ahler $f$-structure
on a $6$-dimensional manifold there exists at least one
corresponding Einstein metric.

{\bf Remark 4.} The invariant $f$-structures $f_4,f_5,f_6$ on the
flag manifold $M=SU(3)/T_{max}$ cannot be canonical $f$-structures
for all homogeneous $\Phi$-spaces of orders $4$ and $5$ of the Lie
group $SU(3)$. It evidently follows by comparing the results in
Theorem \ref{t8}, Theorem \ref{t9}, and Proposition \ref{p3}.

\subsection{\bf Hermitian $f$-structures}

First let us calculate the composition tensor $T$ (see formula
(\ref{eqT})) for arbitrary invariant $f$-structure on
$(M=SU(3)/T_{max},g=(1,t,s))$. Combining (\ref{123}) and
(\ref{eqf}), we can obtain:
\begin{multline}\label{T}
  T(X,Y)=\frac18 D(\check A,\check B,\check C), \text{where}\\
  \check A=-\zeta_1\zeta_2\zeta_3(1+\zeta_2\zeta_3)((\zeta_1+\zeta_3)(1+s-t)\overline{bc_1}+
  (\zeta_1+\zeta_2)(s-t-1)\overline{b_1c}),\\
  \check B=-\zeta_1\zeta_2\zeta_3(1+\zeta_1\zeta_3)((\zeta_2+\zeta_1)(1+\frac{1-s}{t})\overline{ca_1}+
  (\zeta_2+\zeta_3)(\frac{1-s}{t}-1)\overline{c_1a}),\\
  \check C=-\zeta_1\zeta_2\zeta_3(1+\zeta_1\zeta_2)((\zeta_3+\zeta_2)(\frac{t-1}{s}+1)\overline{ab_1}+
  (\zeta_3+\zeta_1)(\frac{t-1}{s}-1)\overline{a_1b}).
\end{multline}
We recall that the defining property for a Hermitian $f$-structure
is the condition $T(X,Y)=0$. Now from (\ref{T}), we get the
following result:

\begin{propos}\label{p4}
The invariant $f$-structures $J_2,J_3,J_4$ and $f_1,\dots,f_9$ are
Hermitian $f$-structures with respect to all invariant Riemannian
metrics $g=(1,t,s)$,\newline$t,s>0$ on the flag manifold
$M=SU(3)/T_{max}$.
\end{propos}

Notice that the almost complex structure $J_1=(1,1,1)$ is
non-integrable. It agrees with the fact that $J_1$ is not a
Hermitian $f$-structure for each Riemannian metric. While we
stress that all $f$-structures $f_1,\dots,f_9$ of rank $4$ and $2$
are non-integrable, but they are Hermitian $f$-structures.

\subsection{\bf G$_1$f-structures}

Finally, we consider the condition $T(X,X)=0$, which is the
defining property for $G_1f$-structures. Using (\ref{T}) and
taking into account Propositions \ref{p3} and \ref{p4}, we get

\begin{propos}\label{p5}
The flag manifold $M=SU(3)/T_{max}$ does not admit invariant
strictly $G_1f$-structures (i.e. neither $NKf$-structures nor
$Hf$-structures). In particular, there are no invariant strictly
$G_1$-structures $J$ (i.e. neither nearly K\"ahler nor Hermitian)
on $M$.
\end{propos}

\bigskip

Vitaly V. Balashchenko \\
Faculty of Mathematics and Mechanics  \\
Belarusian State University    \\
F.Scorina av.~4 \\
Minsk 220050, BELARUS      \\
{\it E-mail}: balashchenko@bsu.by \\
\hspace*{1.4cm} vitbal@tut.by

\end{document}